\documentclass[1p]{elsarticle}
\journal{Statistics and Probability Letters}
\usepackage{natbib}
\usepackage{enumerate}
\usepackage{amsmath,amsthm}
\usepackage{amssymb, latexsym}
\usepackage[mathscr]{euscript}
\usepackage[toc, page]{appendix}
\usepackage{hyperref}
\usepackage{color}
\usepackage{array}
\usepackage{flafter}
\usepackage{layout}
\usepackage{graphicx}
\usepackage{epstopdf}
\usepackage[ansinew]{inputenc}
\usepackage{xfrac}
\usepackage{algorithm}
\theoremstyle{plain}
  \newtheorem{teo}{Theorem}

\theoremstyle{definition}

\theoremstyle{remark}
  \newtheorem{remark}{Remark}

\begin{document}

\author[1]{Daniel Andr\'es D\'{\i}az--Pach\'on\corref{cor1}} \ead{Ddiaz3@miami.edu}
\author[2]{Juan Pablo S\'aenz} \ead{j.saenz4@umiami.edu}
\author[1]{J. Sunil Rao} \ead{JRao@biostat.med.miami.edu}
\cortext[cor1]{Corresponding author}
\address[1]{Division of Biostatistics - University of Miami, Don Soffer Clinical Research Center, 1120 NW 14th St, Miami FL, 33136}
\address[2]{Department of Industrial Engineering - University of Miami, McArthur Engineering Building, 1251 Memorial Drive, Coral Gables FL, 33146}

\begin{abstract}
We develop hypothesis testing for active information ---the averaged quantity in the Kullback-Liebler divergence. To our knowledge, this is the first paper to derive \textit{exact} probabilities of type-I errors for hypothesis testing in the area.
\end{abstract}
\title{Hypothesis testing with active information}
\maketitle


\section{Introduction}

The No Free Lunch Theorems (NFLT)  \cite{WolpertMacready1995, WolpertMacready1997}, introduced by Wolpert and Macready, proved that no search behaves better on average than a blind search. For instance, selecting a point according to a standard normal r.v.\ restricted to the interval $[-x, x]$ might be good at finding a target close to 0, but it won't work as well towards the limits of the interval; other r.v.'s will do better in the latter case.

However, in applications, searches like evolutionary algorithms usually do better than blind chance. Wolpert and Macready attributed it to the incorporation of ``problem-specific knowledge into the behavior of the [optimization or search] algorithm.'' Active information (actinfo) was thus introduced by Dembski and Marks in order to measure this amount of knowledge infused by the programmer to reach a given target \cite{DembskiMarks2009a, DembskiMarks2009b}. Actinfo is obtained by measuring how much information is infused in an algorithm in order to reach a given target and then subtracting from it the information of reaching that same target by a uniformly-at-random search.

Recently, Monta\~nez proposed a model unifying different versions of complex specified information, setting in turn the stage to use actinfo in hypothesis testing \cite{Montanez2018}. In this paper we build on Monta\~nez's idea of hypothesis testing. Thus, the main goal of this article is to set a general framework for hypothesis testing. Accordingly, we begin by doing a basic review of active information and the framework of hypothesis testing developed by Monta\~nez.

\subsection{Active information}

Let's consider a search space $\Omega$ and a target $T\subset \Omega$ such that $|T| \ll |\Omega|$, where $|\cdot|$ stands for cardinality. The NFLT show that no search behaves better on average than $\textbf U(\cdot)$, where $\textbf U$ is uniform over $\Omega$. (This implies that $|\Omega| < \infty$, since in any other case the uniform distribution does not exist.) Let's define the endogenous information as $-\log \textbf U(T)$; we will denote it by $I_0$. Defined this way, $I_0$ measures the inherent difficulty of reaching  $T$.

Of course, alternative searches can  ---and must--- be developed once relevant additional knowledge is acquired. Partial or total knowledge on the position of the target or the space structure might alter the chances of reaching $T$. Such knowledge assigns a new probability $p$ of reaching $T$. We now call $-\log p$ the exogenous information, and we will denote it by $I_1$.

The difference $I_+ := I_0 - I_1 = \log( p/\textbf U)$ is called active information. It measures the amount of information added to the search by the programmer with respect to the one provided by the blind search. 


\subsection{Hypothesis testing with active information}

In order to explain Monta\~nez's ideas in \cite{Montanez2018}, we begin with a null hypothesis that is rejected whenever its $p$-value, $p_\text{val}$, satisfies that $p_\text{val}< \alpha$ for $\alpha \in (0,1)$. Or, equivalently, we reject the null hypothesis when
\begin{align}\label{transfo}
	\log \frac{p_\text{val}}{\alpha} < 0.
\end{align}
Since, under the null hypothesis, $p$-values are uniformly distributed in $(0,1)$ when the test statistic defining the $p$-value is continuous, 
\begin{align}
	\textbf P[p_\text{val} < \alpha] \leq \alpha.
\end{align}
From this, we obtain directly that
\begin{align}
	\textbf P \left[  \log \frac{p_\text{val}}{\alpha} < x \right] \leq \alpha \exp(x).
\end{align}

This last equation corresponds to Theorem 1 of \cite{Montanez2018}. 
The following Theorem can be easily understood from the previous discussion:
\begin{teo}[Theorem 2 of \cite{Montanez2018}: Conservation of canonical specified complexity]
	Let $p(x)$ be any probability measure on space $\mathcal X$, and let $v: \mathcal X \rightarrow \mathbb R^+$ be an integrable function, such that $v(\mathcal X) \leq r$ for a 		constant  $r \in \mathbb R^+$. For $X \sim p$,
	\begin{align}\label{Monta}
		\textbf P \left[ -\log r \frac{p(X)}{v(X)} \geq x\right] \leq \exp(-x).
	\end{align}
\end{teo}

Notice that up to this point we have not specified the base of the logarithm, since in principle it could be given in different units ---bits, trits, dits, nats, etc. Throughout this article, unless we explicitly specify the base, we are not assuming any particular unit.  


\section{Coin tossing}


Imagine that we are tossing a coin. Against the backdrop of a search, $p$ can be thought as the probability of reaching the target ``heads'' under a new search strategy, and $1/2$ as the default probability of reaching ``heads'' in the absence of further knowledge. For these reasons, appealing to the definitions in the previous section, we call them respectively exogenous probability and endogenous probability. When these probabilities are close, $\log p/(1/2) \approx 0$. Of course, we know that if $p$ is too far removed from 1/2, say $p=1$, then the distribution is not uniform. But how close is close enough so that it is not justified to reject the hypothesis that the probability $p$ is effectively the same as that of a uniform?  Active information allows us to link information to rejection regions specifying how many information units we need in an event in order to reject the null hypothesis to a given $\alpha$-level:
\begin{align}\label{activeinfo}
	\textbf H_0: 0 \approx I_+.
\end{align}




\subsection{One-sided test}

In order to be able to use $I_+$ as a test statistic, it has to be random. So let's suppose $p$ is distributed as a continuous uniform random variable  in $(0,1]$. (This is equivalent to set ourselves in a Bayesian framework in which we have a Bernoulli distribution with random parameter $p$ such that the parameter is uniformly distributed in $(0,1]$.) In bits, this distribution is given by: 
\begin{align}\label{UniBin}
	\textbf P \left[I_+ \leq b\right] = \textbf P\left[\log_2 \frac{p}{1/2} \leq b\right] 
		&= 
		\begin{cases}
			2^{b-1}	, \ \ \ \ \ \ \ \ \text{if } b \in (-\infty,1];\\
			1, \ \ \ \ \ \ \ \ \ \ \ \ \text{if } b > 1.
		\end{cases}
\end{align}

This enables us to find the probability of type-one errors exactly:
\begin{align}\label{alpha2bin}
	\textbf P\left[ I_+ > b\right] < \alpha &\Leftrightarrow b > 1 + \log_2 (1-\alpha),
\end{align}
for $0 < \alpha < 1$. Table 1 in the supplementary material gives $\alpha$-levels from (\ref{alpha2bin}).


\subsection{Two-sided test}

Since $I_+$ can also be negative (when $p <1/2)$, it is needed to consider a deviation from the two sides, which in nats becomes:
\begin{align}\label{2nats}
	\textbf P \left[ \left| I_+ \right| \leq n \right] &= 
		 \begin{cases}
				\sinh n, \ \ \ \ \ \ \ \ \ \ \ \ \ \ \ \ n \leq \ln 2;\\
				1 - \frac{\cosh n - \sinh n}{2}, \ \ \ n > \ln 2.
			\end{cases}.
\end{align}

\begin{remark}
	The interesting result obtained in (\ref{2nats}) in terms of the hyperbolic functions explains why we took the base $e$. Notice then that the density function of the r.v. $ Y = |I_+|$ is given by 
	\begin{align}
		f_Y(n) = \begin{cases}
				 \cosh n \ \ \ \ \ \ \ \ \ \ \ \ \  n \leq \ln 2;\\
				 \frac{\cosh n - \sinh n}{2} \ \ \ \ \ \ n > \ln 2.
				\end{cases}
	\end{align}
	In other words, the density of $|I_+|$ can be seen as a catenary function $y = \cosh n$, and its distribution as the arc length of the catenary from 0 to $n$, provided $n \leq \ln 2$.
\end{remark}

Going to nats provides an easy way to find an $\alpha$-value, since the inverses of hyperbolic functions are well-known. We obtain from (\ref{2nats}) that 
\begin{align}
	\textbf P \left[ \left| I_+ \right| > n \right] &=
		\begin{cases}
			1 - \sinh n, \ \ \ \ \ \ n \leq \ln 2;\\
			\frac{e^{-n}}{2}, \ \ \ \ \ \ \ \ \ \ \ \ \ \ n > \ln 2.
		\end{cases}
\end{align}

Thus, for $\textbf P \left[ |I_+| > n\right] < \alpha$, we obtain that
\begin{align}\label{alpha2nats}
	n > \begin{cases}
			\ln \left(1 - \alpha + \sqrt{ (1-\alpha)^2 + 1}\right), \ \ \ \ \ \ \ \ \alpha \leq 1/4; \\
			-\ln (2\alpha), \ \ \ \ \ \ \ \ \ \ \ \ \ \ \ \ \ \ \ \ \ \ \ \ \ \ \ \ \ \ \ \ \ \ \alpha > 1/4.
		\end{cases}
\end{align}

Table 2 in the supplementary material gives some $\alpha$-levels from (\ref{alpha2nats}). 


\section{General discrete uniform $U(N)$}

We have detailed the analysis for a discrete uniform r.v. with two points. But this can be easily generalized to any $N\in \mathbb N$. That is, consider $\Omega = \{1,\ldots, N\}$. Let $X\sim U(N)$, and let $Y$ be an arbitrary r.v. in the same space.

If for a given singleton the probability under $Y$ is $p$, the active information becomes
\begin{equation}
	I_+ = \log_N \frac{p}{1/N} = 1 + \log_N (p).
\end{equation}

Here we consider the base of the logarithm to be $N$, so that we measure our information in ``$N$-its''. This simplifies the notation and can be easily converted to bits by means of the equation $N^x = 2^b$. Assuming $p\sim \mathcal U((0,1])$, the distribution of $I_+$ becomes
\begin{align}\label{UniMulti}
	\textbf P \left[I_+ \leq x \right] &= 
					\begin{cases}
							 N^{x-1}, \ \ \ \ \ \ \ \ x \in (-\infty,1];\\
							 1, \ \ \ \ \ \ \ \ \ \ \ \ \ \  x > 1;
					\end{cases}
\end{align}
From (\ref{UniMulti}), when $ x < 1$,
\begin{align}\label{UniMulti2}
	\textbf P\left[ I_+ > x\right] < \alpha &\Leftrightarrow 1 + \log_N (1-\alpha) < x, 
\end{align}
and $\textbf P \left[ I_+ > x \right] = 0$ whenever $x > 1$. 

On the other hand, in the two-sided case, considered again in nats:
\begin{align}
	\textbf P \left[ \left | I+ \right | \leq n \right] 
		&= \begin{cases}
			\frac{2}{N}\sinh n, \ \ \ \ \ \ \ \ \ \ \ \ \ \ \ \ \ n \leq \ln N;\\
			1- \frac{\cosh n - \sinh n}{N}, \ \ \ \ \ \ \ \ n > \ln N.
		      \end{cases}	
\end{align}
So that for $\textbf P \left[ \left | I+ \right | > n \right] < \alpha$, we obtain
\begin{align*}
	n  > \begin{cases}
		\ln \left( \frac{(1-\alpha)N}{2} + \sqrt{\left (\frac{(1-\alpha)N}{2} \right)^2 + 1}\right), \ \ \ \ \ \ \ \alpha \leq \frac{1}{N};\\
		 -\ln (\alpha N), \ \ \  \ \ \ \ \ \ \ \ \ \ \ \ \ \ \ \ \ \ \ \ \ \ \ \ \ \ \ \ \ \ \ \ \ \ \ \  \alpha > \frac{1}{N}.
		\end{cases}
\end{align*}

\section {General prior $F$}
Our guiding principle has been to set $p$ as a uniform r.v. But what if we want to use a different prior for $p$? For instance, taking Jeffrey's prior, we could consider $p \sim \beta(\frac{1}{2},\frac{1}{2})$. Then, provided that we know the prior distribution for $p$, the actual distribution of $I_+$ can be derived directly. 
%
For instance, for a discrete space with $N$ elements, the test statistic has distribution 
\begin{align*}
	\textbf P[I_+ \leq b] &= F\left(\frac{2^b}{N}\right),\\
	\textbf P[I_+ \leq n] &= F\left(\frac{e^n}{N}\right),
\end{align*}
where the first distribution is given in bits and the second in nats, and $F$ is the prior, the distribution of $p$.

Finally, up to this point we have referenced actinfo to a discrete uniform rv with parameter $N$, since it has maxent over all finite spaces of size $N$, making it a natural choice to compare against \cite{DiazMarks2020a}. In  Section 1 of the supplementary material, we present a general version that removes this condition.

\section{Discussion}

One of the main strengths of actinfo in comparison with other strategies is that it provides an objective way to measure differences: although a quotient of probabilities can be somewhat obscure to interpret, the logarithm of that quotient has a more straight-forward interpretation in terms of bits (or any other information unit) added. Another advantage is that the quotient of the probabilities of a given event of two dimensioned r.v.'s is dimensionless. This seems to support the idea that actinfo is a more basic unit than others in statistics and information theory. In fact, this is highlighted by the fact that the Kullback-Liebler distance ---the average of the actinfo--- is invariant under parameter transformation (see \cite{Kullback1959}, p. 18-22), a fact of relevance both in Bayesian and frequentist statistics.

We have constructed test statistics building from the most simple case of a space with two singletons to the more general situations. Two differences are worth mentioning with respect to Montanez's results in \cite{Montanez2018}. First, Montanez developed some tables similar to ours. However, since he considered $-\log (p/v)$ instead of $\log (p/v)$, the values he obtained differ from ours. It seems to us more intuitive to consider the specification function $v$ in the denominator. In fact, when the specification is the endogenous distribution and the alternate distribution (search) is given by $p$, it is natural to consider this order, since it coincides with what was written in Section 1.1. Second, knowing specifically the distribution of $p$ enables us to find the exact probabilities of the false positives, which is more potent than the inequality in Theorem 1.



Actinfo is at the core of the algorithm called AIMH (active information mode hunting). This algorithm is more efficient to find modes in large dimensions than its competitors, as illustrated in \cite{DiazEtAl2019} (for other models of bump hunting see e.g., \cite{DazardRao2010, DiazDazardRao2017, FriedmanFisher1999}). However, other applications are possible; for instance, actinfo is able to compare two different learning strategies. One of these situations can be envisioned when we want to compare an unsupervised learner to a supervised one. Indeed, in this case the unsupervised learner can be taken as the one to which the actinfo is referenced, and the supervised case might be the one for which we are measuring how much information is added in order to reach the target. We can expect big additions of information of supervised learners that perform better than unsupervised ones. In fact, Section 2 of the Supplementary material ensures that we can compare any two strategies. 

\bibliographystyle{plainnat}
\bibliography{/Users/daniela.diaz/Documents/Research/daangapaBibliography.bib}

\end{document}